\documentclass[11pt]{article}

\usepackage{amsmath,amsthm,amsfonts,bbold}

\usepackage{amssymb}

\usepackage{epsfig}

\usepackage{rotating}

\usepackage{subfigure}

\begin{document}

\title{Path sampling for particle filters with application to multi-target tracking}
\author{Vasileios Maroulas \\
Department of Mathematics \\
University of Tennessee \\
Knoxville, TN 37996 \\
and \\
 Panos Stinis \\
Department of Mathematics \\
University of Minnesota \\
Minneapolis, MN 55455}
%\address{Department of Mathematics,
%    University of California
%and Lawrence Berkeley National Laboratory,
%    Berkeley, CA 94720}
%\email{stinis@math.lbl.gov}
\date {}
%\thanks{This work was supported in part by the Director, Office of Science, Advanced Scientific Computing Research, U.S. Department of Energy under Contract No.
%DE-AC02-05CH11231.}

\maketitle

\begin{abstract}
In recent work \cite{MS10}, we have presented a novel approach for improving particle filters for multi-target tracking. The suggested approach was based on drift homotopy for stochastic differential equations. Drift homotopy was used to design a Markov Chain Monte Carlo step which is appended to the particle filter and aims to bring the particle filter samples closer to the observations. In the current work, we present an alternative way to append a Markov Chain Monte Carlo step to a particle filter to bring the particle filter samples closer to the observations. Both current and previous approaches stem from the general formulation of the filtering problem. We have used the currently proposed approach on the problem of multi-target tracking for both linear and nonlinear observation models. The numerical results show that the suggested approach can improve significantly the performance of a particle filter.
\end{abstract}

%\maketitle

\section*{Introduction}
Multi-target tracking is a central and difficult problem arising
in many scientific and engineering applications including 
radar and signal processing, air traffic control and GPS
navigation \cite{Mah}. The tracking problem consists of computing the 
best estimate of the targets' trajectories based on noisy measurements (observations). 
Several strategies have been developed for addressing the
multi-target tracking problem, see e.g. \cite{bar, FBS,doucet,godsill,gordon,LiCh,Mah,Mah3,MaMa,vo}. 

As in our recent work \cite{MS10}, in this paper we also focus on particle filter techniques \cite{doucet,godsill}.
The popularity of the particle filter method has increased due to its flexibility to handle cases where 
the dynamic and observation models 
are non-linear and/or non-Gaussian. The particle filter
approach is an importance sampling method which approximates the 
target distribution by a discrete set of weighted samples (particles). The weights of the samples are updated when observations become available in order to
incorporate information from the observations.

Despite the particle filter's flexibility, it is often found in practice that most samples will have a
negligible weight with respect to the observation, in other words
their corresponding contribution to the target distribution will
be negligible. Therefore, one may resample the weights to create
more copies of the samples with significant weights
\cite{gordon}. However, even with the resampling step, the
particle filter might still need a lot of samples in order to
approximate accurately the target distribution. Typically, a few 
samples dominate the weight distribution, while the rest of the samples
are in statistically insignificant regions. Thus, some authors (see e.g. \cite{gilks,W09}) have suggested the 
use of an extra step, after the resampling step, which can help move more samples in statistically significant regions.

The extra step for the particle filter is a problem of conditional path sampling for stochastic differential equations (SDEs). In \cite{S10}, a new approach to conditional path sampling based on drift homotopy was presented. In that paper, it was also shown how the algorithm can be used to perform the extra step of a particle filter. In \cite{MS10}, we applied the conditional path sampling algorithm from \cite{S10} to perform the extra step of a particle filter for the problem of multi-target tracking. The numerical results in \cite{MS10} suggested that the approach can improve significantly the performance of a particle filter for multi-target tracking. In the current work, we show yet another way of how to perform the extra step of a particle filter. Both the current approach and the one in \cite{MS10} stem from the general formulation of the filtering problem. The details of the currently proposed implementation of the extra step for a single target are given in Section \ref{MCMC_step}  and for multiple targets in Section \ref{MCMC_step_multi}. The relative merits of the proposed approach in this paper and the one proposed in \cite{MS10} are briefly discussed in Section \ref{discussion}. A more detailed comparison will be presented elsewhere.

To address the target-observation association problem we have used a simple Metropolis Monte
Carlo algorithm which first appeared in \cite{MS10}.
This algorithm effects a probabilistic search of the space of
possible associations to find the best target-observation
association. Of course, one can use more sophisticated association
algorithms (see \cite{godsill} and references therein) but the Monte Carlo algorithm performed very well in
the numerical experiments.

The paper is organized as follows. Sections \ref{generic} and \ref{generic_multi}
provide a brief presentation of particle filters for single and multiple targets (more details
can be found in \cite{gordon,doucet,liu,godsill}), which will serve to
highlight the versatility and drawbacks of this popular filtering
method. Sections \ref{MCMC_step} and \ref{MCMC_step_multi}  demonstrate how one can use an
extra step to improve the performance of particle filters for single and multiple targets. 
%Section \ref{association} describes the
%Monte Carlo sampling algorithm for computing the best association
%between observations and targets for the case of multiple targets.
Section \ref{numerical} presents numerical
results for multi-target tracking for the cases of linear and nonlinear observation models. Finally, Section \ref{discussion} contains 
a discussion of the results as well as directions for future work.

\section{Particle filtering}\label{particle_filtering}
Particle filters are a special case of sequential importance
sampling methods. In Sections \ref{generic} and \ref{generic_multi} we discuss the generic particle 
filter for a single and multiple targets respectively. In Sections \ref{MCMC_step} and \ref{MCMC_step_multi} 
we discuss the addition of an extra step to the generic particle filter for the cases of a single and multiple targets respectively.

\subsection{Generic particle filter for a single target}\label{generic}
Suppose that we are given an SDE system and that we also have
access to noisy observations $Z_{T_1},\ldots,Z_{T_K}$ of the state
of the system at specified instants $T_1,\ldots,T_K.$ The
observations are functions of the state of the system, say given
by $Z_{T_k}=G(X_{T_k},\xi_k),$ where $\xi_k, k=1,\ldots,K$ are
mutually independent random variables. For simplicity, let us
assume that the distribution of the observations admits a density
$g(X_{T_k},Z_{T_k}),$ i.e., $p(Z_{T_k}|X_{T_k} ) \propto
g(X_{T_k},Z_{T_k}).$

The filtering problem consists of computing estimates of the
conditional expectation $E[f(X_{T_k})| \{Z_{T_j}\}^{k}_{j=1}],$
i.e., the conditional expectation of the state of the system given
the (noisy) observations. Equivalently, we are looking to compute
the conditional density of the state of the system given the
observations $p(X_{T_k}|\{Z_{T_j}\}^{k}_{j=1}).$ There are several
ways to compute this conditional density and the associated
conditional expectation but for practical applications they are
rather expensive.

Particle filters fall in the category of importance sampling
methods. Because computing averages with respect to the
conditional density involves the sampling of the conditional
density which can be difficult, importance sampling methods
proceed by sampling a reference density
$q(X_{T_k}|\{Z_{T_j}\}^{k}_{j=1})$ which can be easily sampled and
then compute the weighted sample mean
$$E[f(X_{T_k})| \{Z_{T_j}\}^{k}_{j=1}] \approx \frac{1}{N} \sum_{n=1}^N
f(X^n_{T_k})\frac{p(X^n_{T_k}|\{Z_{T_j}\}^{k}_{j=1})}{q(X^n_{T_k}|\{Z_{T_j}\}^{k}_{j=1})}$$
or the related estimate
\begin{equation}\label{importance}
E[f(X_{T_k})| \{Z_{T_j}\}^{k}_{j=1}] \approx \frac{\sum_{n=1}^N
f(X^n_{T_k})
\frac{p(X^n_{T_k}|\{Z_{T_j}\}^{k}_{j=1})}{q(X^n_{T_k}|\{Z_{T_j}\}^{k}_{j=1})}}{\sum_{n=1}^N
\frac{p(X^n_{T_k}|\{Z_{T_j}\}^{k}_{j=1})}{q(X^n_{T_k}|\{Z_{T_j}\}^{k}_{j=1})}},
\end{equation}
where $N$ has been replaced by the approximation
$$N \approx \sum_{n=1}^N \frac{p(X^n_{T_k}|\{Z_{T_j}\}^{k}_{j=1})}{q(X^n_{T_k}|\{Z_{T_j}\}^{k}_{j=1})}.$$
Particle filtering is a recursive implementation of the importance
sampling approach. It is based on the recursion
\begin{align}
p(X_{T_k}|\{Z_{T_j}\}^{k}_{j=1}) & \propto g(X_{T_k},Z_{T_k}) p(X_{T_k}|\{Z_{T_j}\}^{k-1}_{j=1}),
\label{correct}\\  
\text{where} \;\; p(X_{T_k}|\{Z_{T_j}\}^{k-1}_{j=1}) & =  \int p(X_{T_k}|
X_{T_{k-1}})p(X_{T_{k-1}}|\{Z_{T_j}\}^{k-1}_{j=1}) dX_{T_{k-1}}.
\label{update}
\end{align}
If we set
$$q(X_{T_k}|\{Z_{T_j}\}^{k}_{j=1})=p(X_{T_k}|\{Z_{T_j}\}^{k-1}_{j=1}),$$
then from \eqref{correct} we get
$$\frac{p(X_{T_k}|\{Z_{T_j}\}^{k}_{j=1})}{q(X_{T_k}|\{Z_{T_j}\}^{k}_{j=1})} \propto g(X_{T_k},Z_{T_k}).$$
The approximation in expression \eqref{importance} becomes
\begin{equation}\label{particle}
E[f(X_{T_i})| \{Z_{T_j}\}^{k}_{j=1}] \approx \frac{\sum_{n=1}^N
f(X^n_{T_k})g(X^n_{T_k},Z_{T_k})}{\sum_{n=1}^N
g(X^n_{T_k},Z_{T_k})}
\end{equation}
From \eqref{particle} we see that if we can construct samples from
the predictive distribution $p(X_{T_k}|\{Z_{T_j}\}^{k-1}_{j=1})$
then we can define the (normalized) weights $W^n_{T_k}=
\frac{g(X^n_{T_k},Z_{T_k})}{\sum_{n=1}^N g(X^n_{T_k},Z_{T_k})},$
use them to weigh the samples and the weighted samples will be
distributed according to the posterior distribution
$p(X_{T_k}|\{Z_{T_j}\}^{k}_{j=1}).$

In many applications, most samples will have a negligible weight
with respect to the observation, so carrying them along does not
contribute significantly to the conditional expectation estimate
(this is the problem of degeneracy \cite{liu}). To create larger
diversity one can resample the weights to create more copies of
the samples with significant weights. The particle filter with
resampling is summarized in the following algorithm due to Gordon
{\it et al.} \cite{gordon}.

\vskip14pt
{\bf Particle filter for a single target}
\begin{enumerate}
\item Begin with $N$ unweighted samples $X^n_{T_{k-1}}$ from
$p(X_{T_{k-1}}|\{Z_{T_j}\}^{k-1}_{j=1}).$ \item {\bf Prediction}:
Generate $N$ samples $X'^n_{T_k}$ from $ p(X_{T_k}| X_{T_{k-1}}).$
\item {\bf Update}: Evaluate the weights $$W^n_{T_k}=
\frac{g(X'^n_{T_k},Z_{T_k})}{\sum_{n=1}^N
g(X'^n_{T_k},Z_{T_k})}.$$ 
\item {\bf Resampling}: Generate $N$
independent uniform random variables $\{\theta^n\}_{n=1}^N$ in
$(0,1).$ For $n=1,\ldots,N$ let $X^n_{T_k}=X'^j_{T_k} $where $$
\sum_{l=1}^{j-1}W^l_{T_k} \leq \theta^j < \sum_{l=1}^{j}W^l_{T_k}
$$ 
where $j$ can range from $1$ to $N.$
\item Set $k=k+1$ and proceed to Step 1.
\end{enumerate}

The particle filter algorithm is easy to implement and adapt for
different problems since the only part of the algorithm that
depends on the specific dynamics of the problem is the prediction
step. This has led to the particle filter algorithm's increased
popularity \cite{doucet}. However, even with the resampling step,
the particle filter can still need a lot of samples in order to
describe accurately the conditional density
$p(X_{T_k}|\{Z_{T_j}\}^{k}_{j=1}).$ Snyder {\it et al.}
\cite{snyder} have shown how the particle filter can fail in
simple high dimensional problems because one sample dominates the
weight distribution. The rest of the samples are not in
statistically significant regions. Even worse, as we will show in
the numerical results section, there are simple examples where not
even one sample is in a statistically significant region. In the
next subsection we present how an extra step can be used to
push samples closer to statistically significant regions.

\subsection{Generic particle filter for multiple targets}\label{generic_multi}
Suppose that we have $\lambda=1,\ldots,\Lambda$ targets. Also, for notational simplicity, assume that 
the $\lambda$th target comes from the $\lambda$th observation. Even when this is not the case, we can 
relabel the observations to satisfy this assumption. 
The targets are assumed to evolve independently so that the observation weight of a sample of the vector of targets is the product of the individual 
observation weights of the targets \cite{godsill}. The same is true for the transition density of the vector of targets between 
observations. We denote the vector of targets at observation $T_k$ by
$$ X_{T_k}=(X_{1,T_k},\ldots,X_{\Lambda,T_k})$$
and the observation vector at $T_k$ by
$$ Z_{T_k}=(Z_{1,T_k},\ldots,Z_{\Lambda,T_k}).$$ Also, we can have different observation weight densities $g_{\lambda}, \; \lambda=1,\ldots,\Lambda$ for different targets. However, in the numerical examples we have chosen the same observation weight density for all targets.  

Following \cite{godsill} we can write the particle filter for the case of multiple targets as
\vskip14pt
{\bf Particle filter for multiple targets}
\begin{enumerate}
\item Begin with $N$ unweighted samples $X^n_{T_{k-1}}$ from
$p(X_{T_{k-1}}|\{Z_{T_j}\}^{k-1}_{j=1})=\prod_{\lambda=1}^{\Lambda}p(X_{\lambda,T_{k-1}}|\{Z_{\lambda,T_j}\}^{k-1}_{j=1}).$ 
\item {\bf Prediction}:
Generate $N$ samples $X'^n_{T_k}$ from $$ p(X_{T_k}| X_{T_{k-1}})=\prod_{\lambda=1}^{\Lambda}p(X_{\lambda,T_k}| X_{\lambda,T_{k-1}}).$$
\item {\bf Update}: Evaluate the weights $$W^n_{T_k}=
\frac{\prod_{\lambda=1}^{\Lambda} g_{\lambda}({X'}^n_{\lambda,T_k},Z_{\lambda,T_k})}{\sum_{n=1}^N
\prod_{\lambda=1}^{\Lambda} g_{\lambda}({X'}^n_{\lambda,T_k},Z_{\lambda,T_k})}.$$ 
\item {\bf Resampling}: Generate $N$
independent uniform random variables $\{\theta^n\}_{n=1}^N$ in
$(0,1).$ For $n=1,\ldots,N$ let $X^n_{T_k}=X'^j_{T_k} $where $$
\sum_{l=1}^{j-1}W^l_{T_k} \leq \theta^j < \sum_{l=1}^{j}W^l_{T_k}
$$ 
where $j$ can range from $1$ to $N.$
\item Set $k=k+1$ and proceed to Step 1.
\end{enumerate}

\subsection{Particle filter with MCMC step for a single target}\label{MCMC_step}
Several authors (see e.g. \cite{gilks,W09}) have suggested the use
of a MCMC step after the resampling step (Step 4) in order to move
samples away from statistically insignificant regions. There are
many possible ways to append an MCMC step after the resampling
step in order to achieve that objective. The important point is
that the MCMC step must preserve the conditional density
$p(X_{T_k}|\{Z_{T_j}\}^{k}_{j=1}).$ In the current section we show
that the MCMC step constitutes a case of conditional path
sampling. 

We begin by noting that one can use the resampling step (Step 4)
in the particle filter algorithm to create more copies not only of
the good samples according to the observation, but also of the
values (initial conditions) of the samples at the previous
observation. These values are the ones who have evolved into good
samples for the current observation (see more details in
\cite{W09}). The motivation behind producing more copies of the
pairs of initial and final conditions is to use the good initial
conditions as starting points to produce statistically more
significant samples according to the current observation. This
process can be accomplished in two steps. First, Step 4 of the
particle filter algorithm is replaced by

\vskip14pt {\bf Resampling}: Generate $N$ independent uniform
random variables $\{\theta^n\}_{n=1}^N$ in $(0,1).$ For
$n=1,\ldots,N$ let
$(X^n_{T_{k-1}},X^n_{T_k})=(X'^j_{T_{k-1}},X'^j_{T_k}) $where $$
\sum_{l=1}^{j-1}W^l_{T_k} \leq \theta^j < \sum_{l=1}^{j}W^l_{T_k}
$$
Also, through Bayes' rule \cite{W09} one can show that the posterior density $p(X_{T_k}|\{Z_{T_j}\}^{k}_{j=1})$ is 
preserved if one samples from the density $$g(X_{T_k},Z_{T_k}) p(X_{T_k}|X_{T_{k-1}})$$
where $X_{T_{k-1}}$ are given by the modified resampling step. This is a problem of conditional path sampling for (continuous-time or discrete) stochastic systems. The important issue is to perform the necessary sampling efficiently \cite{CT09,W09}. 
  
We are now in a position to present the particle filter with MCMC step algorithm

\vskip14pt
{\bf Particle filter with MCMC step for a single target}
\begin{enumerate}
\item Begin with $N$ unweighted samples $X^n_{T_{k-1}}$ from
$p(X_{T_{k-1}}|\{Z_{T_j}\}^{k-1}_{j=1}).$ 
\item {\bf Prediction}:
Generate $N$ samples $X'^n_{T_k}$ from $ p(X_{T_k}| X_{T_{k-1}}).$
\item {\bf Update}: Evaluate the weights $$W^n_{T_k}=
\frac{g(X'^n_{T_k},Z_{T_k})}{\sum_{n=1}^N
g(X'^n_{T_k},Z_{T_k})}.$$ 
\item {\bf Resampling}: Generate $N$
independent uniform random variables $\{\theta^n\}_{n=1}^N$ in
$(0,1).$ For $n=1,\ldots,N$ let
$(X^n_{T_{k-1}},X^n_{T_k})=(X'^j_{T_{k-1}},X'^j_{T_k}) $ where $$
\sum_{l=1}^{j-1}W^l_{T_k} \leq \theta^j < \sum_{l=1}^{j}W^l_{T_k}
$$
where $j$ can range from $1$ to $N.$
 \item {\bf MCMC step}: For $n=1,\ldots,N$ choose a modified
drift (possibly different for each $n$). Construct a Markov chain
for $Y^{n}_{T_k}$ with stationary distribution
$$ g(Y^n,Z_{T_k}) p(Y^n | X^n_{T_{k-1}}) $$
\item
Set $X^n_{T_k}=Y^{n}_{T_k}.$
\item
Set $k=k+1$ and proceed to Step 1.
\end{enumerate}
Since the samples $X^n_{T_k}=Y^{n,\Lambda}_{T_k}$ are constructed
by starting from different sample paths, they are independent.
Also, note that the samples $X^n_{T_k}$ are unweighted. However, we
can still measure how well these samples approximate the posterior
density by comparing the effective sample sizes of the particle
filter with and without the MCMC step. For a collection of $N$
samples the effective sample size $ess(T_k)$ is defined by
$$ess(T_k) = \frac{N}{1+C_k^2}$$
where
\begin{equation*}
C_k =  \frac{1}{W_k}
\sqrt{\frac{1}{N}\sum_{n=1}^N (g(X^n_{T_k},Z_{T_k}) -W_k)^2} \;\; \text{and} \;\; W_k = \frac{1}{N} \sum_{n=1}^N g(X^n_{T_k},Z_{T_k}).
\end{equation*}
The effective sample size can be interpreted as that the $N$ weighted samples are worth of $ess(T_k) = \frac{N}{1+C_k^2}$ i.i.d. samples drawn from the target density, which in our case is the posterior density. By definition, $ess(T_k) \le N.$ If the samples have uniform weights, then $ess(T_k)=N.$ On the other hand, if all samples but one have zero weights, then $ess(T_k)=1.$

\subsection{Particle filter with MCMC step for multiple targets}\label{MCMC_step_multi}

We discuss now the case of multiple, say $\Lambda,$ targets.
Instead of the observations for a single target now we have a
collection of observations for all the targets
$\{Z_{T_j}\}^{k}_{j=1}=\{(Z^1_{T_j},\ldots,Z^{\lambda}_{T_j})
\}^{k}_{j=1}.$

The particle filter with MCMC step for the case of multiple targets is

\vskip14pt
{\bf Particle filter with MCMC step for multiple targets}
\begin{enumerate}
\item Begin with $N$ unweighted samples $X^n_{T_{k-1}}$ from
$p(X_{T_{k-1}}|\{Z_{T_j}\}^{k-1}_{j=1})=\prod_{\lambda=1}^{\Lambda}p(X_{\lambda,T_{k-1}}|\{Z_{\lambda,T_j}\}^{k-1}_{j=1}).$ 
\item {\bf Prediction}:
Generate $N$ samples $X'^n_{T_k}$ from $$ p(X_{T_k}| X_{T_{k-1}})=\prod_{\lambda=1}^{\Lambda}p(X_{\lambda,T_k}| X_{\lambda,T_{k-1}}).$$
\item {\bf Update}: Evaluate the weights $$W^n_{T_k}=
\frac{\prod_{\lambda=1}^{\Lambda} g_{\lambda}({X'}^n_{\lambda,T_k},Z_{\lambda,T_k})}{\sum_{n=1}^N
\prod_{\lambda=1}^{\Lambda} g_{\lambda}({X'}^n_{\lambda,T_k},Z_{\lambda,T_k})}.$$ 
\item {\bf Resampling}: Generate $N$
independent uniform random variables $\{\theta^n\}_{n=1}^N$ in
$(0,1).$ For $n=1,\ldots,N$ let
$(X^n_{T_{k-1}},X^n_{T_k})=(X'^j_{T_{k-1}},X'^j_{T_k}) $ where $$
\sum_{l=1}^{j-1}W^l_{T_k} \leq \theta^j < \sum_{l=1}^{j}W^l_{T_k}
$$
where $j$ can range from $1$ to $N.$
 \item {\bf MCMC step}: For $n=1,\ldots,N$ and $\lambda=1,\ldots,\Lambda$ choose a modified
drift (possibly different for each $n$ and each $\lambda$). Construct a Markov chain
for $Y^{n}_{T_k}$  with stationary distribution
$$ \prod_{\lambda=1}^{\Lambda} g_{\lambda}(Y^n_{\lambda},Z_{\lambda,T_k}) p_{\lambda}(Y^n_{\lambda}|X^n_{{\lambda},T_{k-1}}). $$
\item
Set $X^n_{T_k}=Y^{n}_{T_k}.$
\item
Set $k=k+1$ and proceed to Step 1.
\end{enumerate}

For a collection of $N$ samples the effective sample size $ess_{\Lambda}(T_k)$ for $\Lambda$ targets is 
$$ess_{\Lambda}(T_k) = \frac{N}{1+C_{\Lambda,k}^2}$$
where
\begin{multline*}
C_{\Lambda,k}=  \frac{1}{W_{\Lambda,k}}
\sqrt{\frac{1}{N}\sum_{n=1}^N (\prod_{\lambda=1}^{\Lambda}g_{\lambda}(X^n_{\lambda,T_k},Z_{\lambda,T_k}) -W_{\Lambda,k})^2}  \\ \text{and} \;\; W_{\Lambda,k} = \frac{1}{N} \sum_{n=1}^N \prod_{\lambda=1}^{\Lambda} g_{\lambda}(X^n_{\lambda,T_k},Z_{\lambda,T_k}).
\end{multline*}

\section{Numerical results}\label{numerical}
We present numerical results for multi-target tracking using the
particle filter with an MCMC step performed by hybrid Monte Carlo. We have synthesized tracks of targets
moving on the $xy$ plane using a $2D$ near constant velocity model
\cite{bar}. At each time $t$ we have a total of $\Lambda_t$ targets and
the evolution of the $k$th target ($\lambda=1,\ldots,\Lambda_t$) is given by
\begin{align}\label{model}
{\bf{x}}_{\lambda,t} &={\bf{A}} {\bf{x}}_{\lambda,t-1} + {\bf{B}}  {\bf{v}}_{\lambda,t} \\
&= [x_{\lambda,t},\dot{x}_{\lambda,t},y_{\lambda,t},\dot{y}_{\lambda,t}]^T, \notag
\end{align}
where $(x_{\lambda,t},\dot{x}_{\lambda,t})$ and $(y_{\lambda,t},\dot{y}_{\lambda,t})$ are
the $xy$ position and velocity of the $k$th target at time $t.$
The matrices ${\bf{A}}$ and ${\bf{B}}$ are given by
\begin{equation}\label{model_matrix1}
{\bf{A}}=
\left[ \begin{array}{cccc}
1 & T & 0 & 0\\
0 & 1 & 0 & 0\\
0 & 0 & 1 & T \\
0 & 0 & 0 & 1
\end{array} \right]
\;\; \text{and} \;\;
{\bf{B}}=
\left[ \begin{array}{cc}
T^2/2 & 0\\
 T & 0\\
0 & T^2/2 \\
0 & T
\end{array} \right],
\end{equation}
where $T$ is the time between observations. For the experiments we
have set $T=1,$ i.e., noisy observations of the model are obtained
at every step of the model \eqref{model}. The model noise
${\bf{v}}_{\lambda,t}$ is a collection of independent Gaussian random
variables with covariance matrix ${\bf{\Sigma}}_{v}$ defined as
 \begin{equation}\label{model_matrix2}
 {\bf{\Sigma}}_{v}=
\left[ \begin{array}{cc}
\sigma_x^2 & 0\\
 0 & \sigma_y^2\\
\end{array} \right].
 \end{equation}
In the experiments we have $\sigma_x^2=\sigma_y^2=1.$ Also, we have considered two possible cases for the observation model, one linear and one nonlinear. Due to the different possible combinations of targets to observations we use a different index $m$ to denote the obsevations. Since we do not assume any clutter we have $m=1,\ldots,\Lambda_t.$ If the $m$th observation ${\bf{z}}_{m,t}$ at time $t$ comes from the $k$th target we have
\begin{equation}\label{linear_observation}
{\bf{z}}_{m,t}=
\left[ \begin{array}{c}
x_{\lambda,t}\\
y_{\lambda,t}\\
\end{array} \right]
+{\bf{w}}_{m,t}
\end{equation}
for the linear observation model and
\begin{equation}\label{nonlinear_observation}
{\bf{z}}_{m,t}=
\left[ \begin{array}{c}
\arctan (\frac{y_{\lambda,t}}{x_{\lambda,t}})\\
(x_{\lambda,t}^2 + y_{\lambda,t}^2)^{1/2}\\
\end{array} \right]
+{\bf{w}}_{m,t}
\end{equation}
for the nonlinear observation model. As is usual in the literature, the nonlinear observation model consists of the bearing $\theta$ and range $r$ of a target. The observation noise ${\bf{w}}_{m,t}$ is white and Gaussian with covariance matrix
 \begin{equation}\label{linear_observation_matrix}
 {\bf{\Sigma}}_{w}=
\left[ \begin{array}{cc}
\sigma_{obs,x}^2 & 0\\
 0 & \sigma_{obs,y}^2\\
\end{array} \right]
 \end{equation}
 for the linear observation model and
  \begin{equation}\label{nonlinear_observation_matrix}
 {\bf{\Sigma}}_{w}=
\left[ \begin{array}{cc}
\sigma_{\theta}^2 & 0\\
 0 & \sigma_{r}^2\\
\end{array} \right]
 \end{equation}
for the nonlinear observation model. For the numerical experiments
with the linear observation model we chose
$\sigma_{obs,x}^2=\sigma_{obs,y}^2=1.$ For the numerical
experiments with the nonlinear observation model we chose
$\sigma_{\theta}^2=10^{-4}$ and $\sigma_{r}^2=1.$ These values
make our example comparable in difficulty to examples appearing in
the literature (see e.g. \cite{godsill,godsill2,vo}).

The synthesized target tracks were created by specifying a certain
scenario, to be detailed below, of surviving, newborn and
disappearing targets. According to this scenario we evolved the
appropriate number of targets according to \eqref{model} and
recorded the state of each target at each step. For the surviving
targets we created an observation by using the state of the target
in  the observation model. Thus, for the linear observation model,
the observations were created directly in $xy$ space by perturbing
the $xy$ position of the target by \eqref{linear_observation}. For
the nonlinear observation model, the observations were created in
bearing and range space $\theta, r$ by using
\eqref{nonlinear_observation}. The perturbed bearing and range
were transformed to $xy$ space by the transformation 
$ x=r \cos \theta,  \; y = r \sin \theta $
to create a position for the
target in $xy$ space.

The newborn targets for the linear model were created in $xy$
space directly by sampling uniformly in $[-100,100].$ Afterwards,
the observations of the newborn targets were constructed by
perturbing the $x,y$ positions using \eqref{linear_observation}.
The newborn targets for the nonlinear model were created in $xy$
space by sampling uniformly in $[-100,100].$ Afterwards, we
transformed the $x,y$ positions to the bearing and range space
$\theta,r$  and perturbed the bearing and range according to
\eqref{nonlinear_observation}. The perturbed bearing and range
were again transformed back to $xy$ space to create the position
of the newborn target. Note that both observation models do {\it
not} involve the velocities. The newborn target velocities were
sampled uniformly in $[-1,1].$
%Finally, in the process of constructing the synthesised tracks, we shuffled the indices of the observation vector at each observation so as to avoid any favorable bias in the target-observation association step.

The  number of targets at each observation instant is: $\Lambda_0=2,$ $\Lambda_1=2$, $\Lambda_2=1,$ $\Lambda_3=2,$ $\Lambda_4=3,$ $\Lambda_5=\Lambda_6=\ldots=\Lambda_{200}=4.$ So, for the majority of the steps we have $4$ targets which makes the problem of tracking rather difficult.

\subsection{MCMC step}

For the $n$th sample, the density we have to sample for the {\it linear} model is
\begin{multline*}
\prod_{\lambda=1}^{\Lambda_{t_k}} g_x(z^n_{1,\lambda,k},Z_{1,\lambda,k})g_y(z^n_{3,\lambda,k},Z_{3,\lambda,k}) p(z^n_{\lambda,k}|z^n_{\lambda,k-1}) \\
\propto \prod_{\lambda=1}^{\Lambda_{t_k}}  \exp \biggl( - \biggl\{  \frac{( Z_{1,\lambda,k}-z^n_{1,\lambda,k-1} - z^n_{2,\lambda,k-1} T -\frac{T^2}{2} v^n_{x,\lambda,k})^2}{2\sigma^2_{obs,x}}  \\
+ \frac{( Z_{3,\lambda,k}-z^n_{3,\lambda,k-1} - z^n_{4,\lambda,k-1} T  -\frac{T^2}{2}  v^n_{y,\lambda,k})^2}{2\sigma^2_{obs,y}} \\
+ \frac{(v^n_{x,\lambda,k})^2}{2\sigma^2_x} + \frac{(v^n_{y,\lambda,k})^2}{2\sigma^2_y} \biggr\}   \biggr).
\end{multline*}

\subsubsection{Hybrid Monte Carlo}
We chose to use Hybrid Monte Carlo to perform the sampling (any other MCMC method can be used) \cite{liu}. 
We present briefly the hybrid Monte Carlo (HMC) formulation that we have used to sample the conditional density for each target. 

Define the potential $V(v^n_{x,k}, v^n_{y,k})$ by
\begin{multline}\label{conditional_linear_potential}
V(v^n_{x,1,k}, v^n_{y,1,k},\ldots,v^n_{x,\Lambda_{t_k},k}, v^n_{y,\Lambda_{t_k},k})= \\
\sum_{\lambda=1}^{\Lambda_{t_k}}  \biggl\{  \frac{( Z_{1,\lambda,k}-z^n_{1,\lambda,k-1} - z^n_{2,\lambda,k-1} T  -\frac{T^2}{2} v^n_{x,\lambda,k})^2}{2\sigma^2_{obs,x}}  \\
+ \frac{( Z_{3,\lambda,k}-z^n_{3,\lambda,k-1} - z^n_{4,\lambda,k-1} T -\frac{T^2}{2}  v^n_{y,\lambda,k})^2}{2\sigma^2_{obs,y}} \\
+ \frac{(v^n_{x,\lambda,k})^2}{2\sigma^2_x} + \frac{(v^n_{y,\lambda,k})^2}{2\sigma^2_y}\biggl\}
\end{multline}
Define the vectors $v^n_{x,k}=[v^n_{x,1,k},\ldots,v^n_{x,\Lambda_{t_k},k}]^T$ and $v^n_{y,k}=[v^n_{y,1,k},\ldots,v^n_{y,\Lambda_{t_k},k}]^T,$ where $T$ is the transpose (not to be confuse with the interval between observations $T$ used before). Consider $v^n_{x,k}, v^n_{y,k}$ as the position variables of a Hamiltonian system. We define the $2\Lambda_{t_k}$-dimensional position vector $q=[q_1,q_2]^T$ with $q_1=v^n_{x,k}$ and $q_2=v^n_{y,k}.$ To each of the position variables we associate a momentum variable and we write the Hamiltonian
$$H (q,p)=V(q)   + \frac{p^T p}{2},$$
where $p=[p_1,p_2]^T$ is the momentum vector. Thus, the momenta variables are Gaussian distributed random variables with mean zero and variance 1. The equations of motion for this Hamiltonian system are given by Hamilton's equations
$$\frac{dq_i}{d\tau}=\frac{\partial H}{\partial p_i} \;\; \text{and} \;\; \frac{dp_i}{d\tau}=-\frac{\partial H}{\partial q_i} \;\; \text{for} \;\; i=1,2.$$
HMC proceeds by assigning initial conditions to the momenta
variables (through sampling  from $\exp(- \frac{p^T p}{2})$),
evolving the Hamiltonian system in fictitious time $\tau$ for a
given number of steps of size $\delta \tau$ and then using the
solution of the system to perform a Metropolis accept/reject step
(more details in \cite{liu}). After the Metropolis step, the
momenta values are discarded. The most popular method for solving
the Hamiltonian system, which is the one we also used, is the
Verlet leapfrog scheme. In our numerical implementation, we did
not attempt to optimize the performance of the HMC algorithm. For
the sampling  we used
$100$ Metropolis accept/reject steps and $1$ HMC step of size
$\delta \tau = 10^{-1}$ to construct a trial path. 

For the nonlinear observation model we can use
the same procedure as in the linear observation model to define a
Hamiltonian system and its associated equations. We omit the
details.

\subsection{Linear observations}\label{Linear}

\begin{figure}
\centering
\epsfig{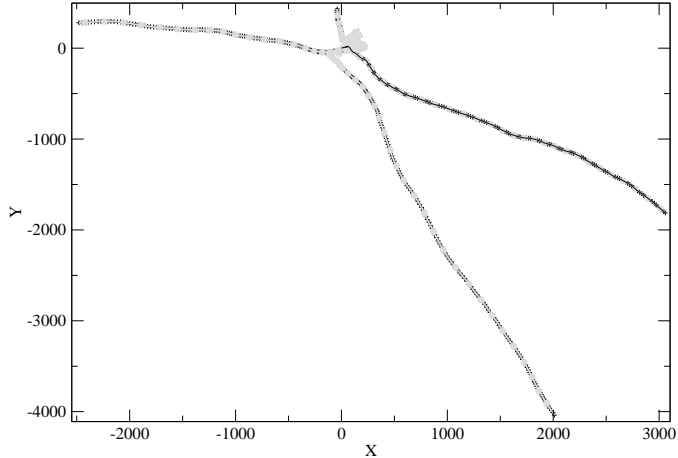}
\caption{Linear observation model. The solid lines denote the true target tracks, the crosses denote the observations and the dots the conditional expectation estimates from the improved particle filter. We have plotted the conditional expectation estimates every 5 observations to avoid cluttering in the figure.}
\label{plot_linear_linear}
\end{figure}

\begin{figure}
\centering
\epsfig{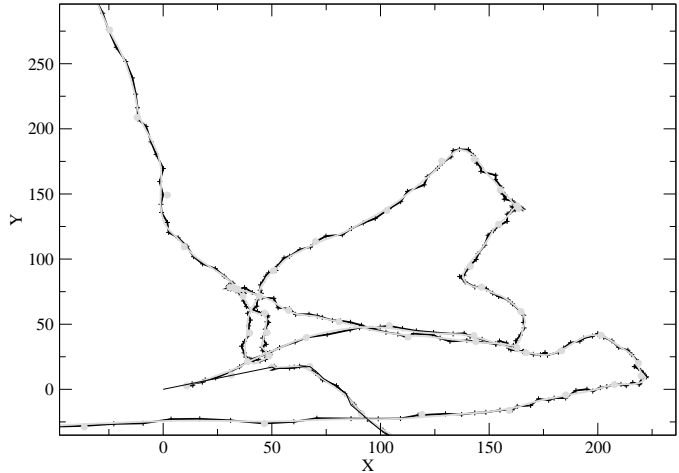}
\caption{Linear observation model. Detail of Figure \ref{plot_linear_linear}. }
\label{plot_linear_linear_focus}
\end{figure}

We start the presentation of our numerical experiments with
results for the linear observation model
\eqref{linear_observation}. Figures \ref{plot_linear_linear} and
\ref{plot_linear_linear_focus} show the evolution in the $xy$
space of the true targets, the observations as well as the
estimates of the improved particle filter. It is obvious from the
figures that the improved particle filter follows accurately the
targets and there is no ambiguity in the identification of the
target tracks. 

\begin{figure}
\centering \epsfig{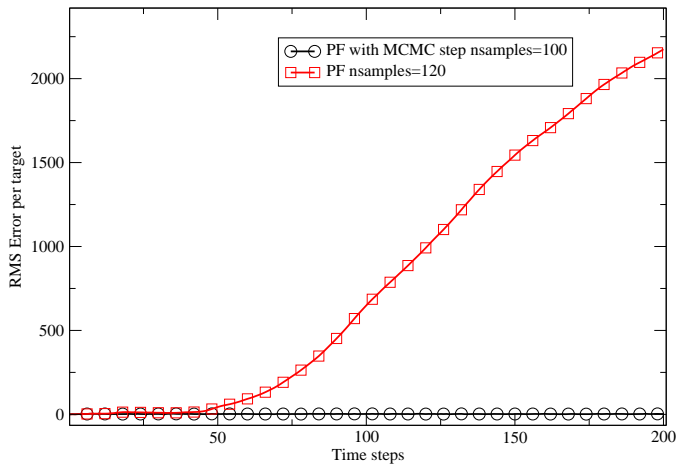}
\caption{Linear observation model. Comparison of RMS error per
target for the improved particle filter and the generic particle
filter.} \label{plot_linear_linear_error}
\end{figure}

The performance of the improved particle filter
with 100 samples is compared to the performance of the generic
particle filter with 120 samples in Figure
\ref{plot_linear_linear_error} by monitoring the evolution in time
of the RMS error per target. The RMS error per target (RMSE) is
defined with reference to the true target tracks by the formula
\begin{equation}\label{rms_error}
RMSE(t) = \sqrt { \frac{1}{K_t}\sum_{k=1}^{K_t}  \|
{\bf{x}}_{k,t} - E[{\bf{x}}_{k,t} | Z_1,\ldots,Z_t ] \|^2    }
\end{equation}
where $\|  \cdot \|$ is the norm of the position and velocity
vector. Note that the state vector norm involves both positions
and velocities even though the observations use information only
from the positions of a target.  ${\bf{x}}_{k,t}$ is the true
state vector for target $k.$ $E[{\bf{x}}_{k,t} |
Z_1,\ldots,Z_t ]$ is the conditional expectation estimate
calculated with the improved or generic particle filter depending on whose 
filter's performance we want to calculate. 

The improved particle filter has a computational overhead of the order of a few
percent compared to the generic particle filter. We have thus used
the generic particle filter with more samples than the improved
particle filter. This additional number of samples more than accounts for the computational
overhead of the improved particle filter. As can be seen in Figure
\ref{plot_linear_linear_error} the generic particle filter's
accuracy deteriorates quickly. On the other hand, the improved
particle filter maintains an $O(1)$ RMS error per target  for the
entire tracking interval. The average value of the RMS error over
the entire time interval of tracking is about 2.5 with standard
deviation of about 0.5. For the generic particle filter, the
average of the RMS error over the time interval of tracking is
about 800 with standard deviation of about 760.

\begin{figure}
\centering \epsfig{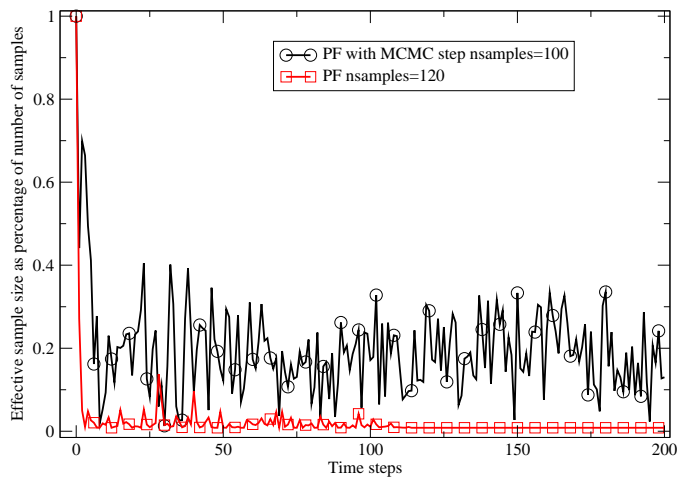}
\caption{Linear observation model. Comparison of effective sample
size for the improved particle filter and the generic particle
filter.} \label{plot_linear_linear_effective}
\end{figure}

Figure \ref{plot_linear_linear_effective} compares the effective
sample size for the generic particle filter and the improved
particle filter. Because the number of samples is different for
the two filters we have plotted the effective sample size as a
percentage of the number of samples. We have to note that, after
about 50 steps, the generic particle filter started producing
observation weights (before the normalization) which were
numerically zero. This makes the normalization impossible. In
order to allow the generic particle filter to continue we chose at
random one of the samples, since all of them are equally bad, and
assigned all the weight to this sample. We did that for all the
steps for which the observation weights were zero before the
normalization. As a result, the effective sample size for the
generic particle filter drops down to 1 sample after about 50
steps. Once the generic particle filter deviates from the true
target tracks there is no mechanism to correct it. Also, we tried
assigning equal weights to all the samples when the observation
weights dropped to zero. This did not improve the generic
particle's performance either. On the other hand, the improved
particle filter maintains an effective sample size which is  about
$25\%$ of the number of samples.

\subsection{Nonlinear observations}\label{Nonlinear}

\begin{figure}
\centering \epsfig{file=Linear_Nonlinear_bw.eps,width=3.5in}
\caption{Noninear observation model. The solid lines denote the
true target tracks, the crosses denote the observations and the
dots the conditional expectation estimates from the improved
particle filter. We have plotted the conditional expectation
estimates every 5 observations to avoid cluttering in the figure.}
\label{plot_linear_nonlinear}
\end{figure}

\begin{figure}
\centering
\epsfig{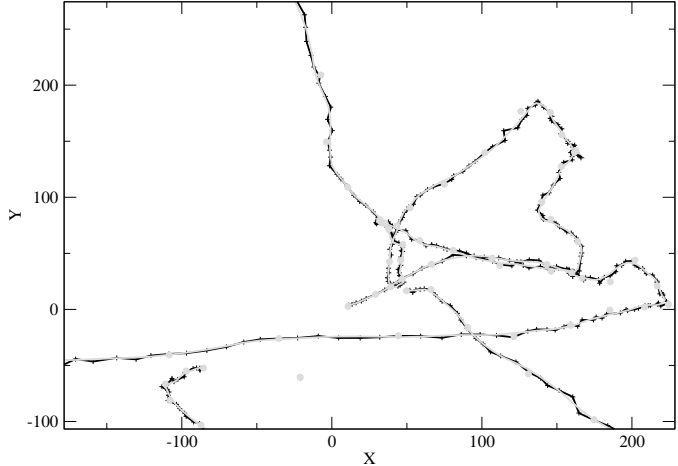}
\caption{Nonlinear observation model. Detail of Figure \ref{plot_linear_nonlinear}. }
\label{plot_linear_nonlinear_focus}
\end{figure}

We continue with results for the nonlinear observation model
\eqref{nonlinear_observation}. Figures \ref{plot_linear_nonlinear}
and \ref{plot_linear_nonlinear_focus} show the evolution in the
$xy$ space of the true targets, the observations as well as the
estimates of the improved particle filter. Again, as in the case
of the linear observation model, the improved particle filter
follows accurately the targets and there is no ambiguity in the
identification of the target tracks.

The case of the nonlinear observation model is much more difficult
than the case of the linear observation model. The reason is that
for the nonlinear observation model, the observation errors,
though constant in bearing and range space, they become position
dependent in $xy$ space. In particular, when $x$ and/or $y$ are
large, the observation errors can become rather large. This is
easy to see by Taylor expanding the nonlinear transformation from
bearing and range space to $xy$ space around the true target
values. Suppose that the true target bearing and range are
$\theta_0, r_0$ and its $xy$ space position is $x_0 = r_0 \cos
\theta_0, y_0= r_0 \sin \theta_0.$ Also, assume that the
observation error in bearing and range space is, respectively,
$\delta \theta$ and $\delta r.$ The $xy$ position of a target that
is perturbed by  $\delta \theta$ and $\delta r$ in bearing and
range space is (to first order)
\begin{align*}
x & =x_0 - y_0 \delta \theta - \delta r \cos \theta_0 \\
y & =y_0 + x_0 \delta \theta  - \delta r \sin \theta_0.
\end{align*}
Thus, the perturbation in $xy$ space can be significant even if
$\delta \theta$ and $ \delta r$ are small. In our example we have
$\sigma_{\theta}=10^{-2}.$ So, when the true target $x$ and $y$
values become of the order of $10^3$ as happens for some of the
targets, the observation value in bearing and range space can be
quite misleading as far as the $xy$ space position of the target
is concerned. As a result, even if one does a good job in
following the observation in bearing and range space, the
conditional expectation estimate of the $xy$ space position can be
inaccurate.

\begin{figure}
\centering \epsfig{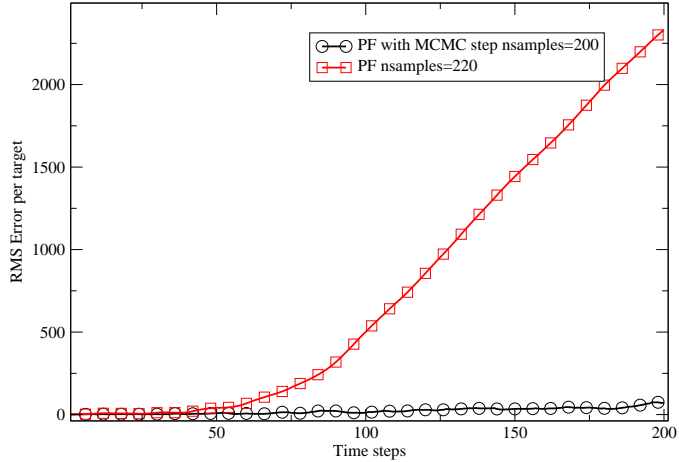}
\caption{Nonlinear observation model. Comparison of RMS error per
target for the improved particle filter and the generic particle
filter.} \label{plot_linear_nonlinear_error}
\end{figure}

With this in mind, we have used 200 samples for the improved
particle filter and 220 samples for the generic particle filter.
Again, the extra samples used for the generic particle filter more
than account for the computational overhead of the improved
particle filter. The performance of the improved particle filter
is compared to the performance of the generic particle filter in
Figure \ref{plot_linear_nonlinear_error} by monitoring the
evolution in time of the RMS error per target. The generic
particle filter's accuracy again deteriorates rather quickly. The
error for the improved particle filter is larger than in the
linear observation model but never exceeds about 80 even after 200
steps when the targets have reached large values of $x$ and/or
$y.$ The average value of the RMS error over the entire time
interval of tracking is about 22 with standard deviation of about
21. For the generic particle filter, the average of the RMS error
over the time interval of tracking is about 760 with standard
deviation of about 770.

\begin{figure}
\centering
\epsfig{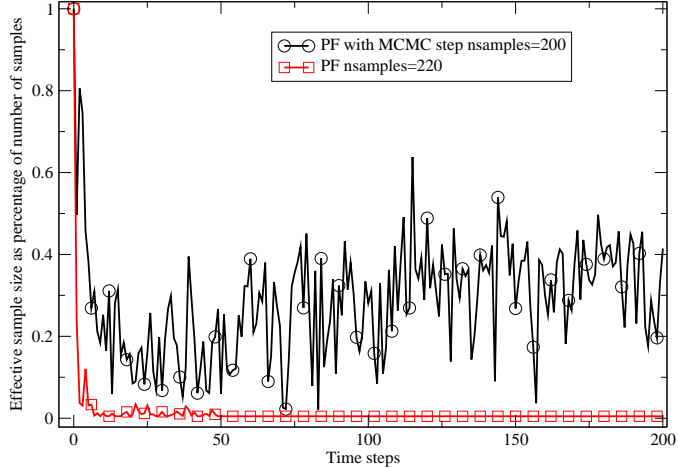}
\caption{Nonlinear observation model. Comparison of effective
sample size for the improved particle filter and the generic
particle filter.} \label{plot_linear_nonlinear_effective}
\end{figure}

Figure \ref{plot_linear_nonlinear_effective} compares the
effective sample size for the generic particle filter and the
improved particle filter. After about 60 steps, the generic
particle filter, started producing observation weights (before the
normalization) which were numerically zero. This makes the
normalization impossible. In order to allow the generic particle
filter to continue we chose at random one of the samples, since
all of them are equally bad, and assigned all the weight to this
sample. We did that for all the steps for which the observation
weights were zero before the normalization. As a result, the
effective sample size for the generic particle filter drops down
to 1 sample after about 60 steps. Once the generic particle filter
deviates from the true target tracks there is no mechanism to
correct it. Also, we tried assigning equal weights to all the
samples when the observation weights dropped to zero. This did not
improve the generic particle's performance either. On the other
hand, the improved particle filter maintains an effective sample
size which is  about $25\%$ of the number of samples.

\section{Discussion}\label{discussion}

We have presented an algorithm for multi-target tracking which  builds on the
existing particle filter methodology for multi-target tracking by
appending an MCMC step after the particle filter resampling step.
The purpose of the addition of the MCMC step is to bring the
samples closer to the observation. Even though the addition of an
MCMC step for a particle filter has been proposed and used before
\cite{gilks}, to the best of our knowledge, the currently proposed implementation of the MCMC step is novel (see also \cite{W09} for a related approach).

We have tested the performance of the algorithm on the problem of
tracking multiple targets evolving under the near constant
velocity model \cite{bar}. We have examined two cases of
observation models: i) a linear observations model involving the
positions of the targets and ii) a nonlinear observation model
involving the bearing and range of the targets. For both cases the
proposed improved particle filter exhibited a significantly better
performance than the generic particle filter. Since the improved
particle filter requires more computations than the generic
particle filter it is bound to be more expensive. However, the
computational overhead of the improved particle filter is rather
small, of the order of a few extra samples worth for the generic
particle filter.

In \cite{MS10} we proposed another way of performing the extra MCMC step of a particle filter. That approach was based on modifying the drift of the dynamic model and then accounting for the modification. In the current work, we use the original drift of the dynamic model without any modification. For the case of multi-target tracking with observations at every time step both algorithms perform equally well. Thus, at first sight it would appear that there is no need for the extra complication of modifying the drift of the dynamic model and then accounting for the modification as was done in \cite{MS10}. However, in cases where there are only sparse observations, the sampling of the conditional density needed for the extra step can be much more difficult (and consequently expensive) for the original dynamic model than for the modified dynamic model. With this in mind, the approach in \cite{MS10} seems to  have wider applicability. A detailed comparison of the algorithm proposed in the current work and the one proposed in \cite{MS10} will be presented in a future publication.

\section*{Acknowledgements}
We are grateful to Prof. J. Weare for many discussions. Also, we would like to thank the Institute for Mathematics and its Applications in the University of Minnesota for its support.

\end{document}